\documentclass[12pt]{article}
\usepackage[latin1]{inputenc}
\usepackage{amsmath, amssymb, amstext, amsfonts, amsbsy, amsthm, fancyhdr, multicol,amscd}
\usepackage{graphicx}
\usepackage{amsfonts, amscd, amssymb}
\usepackage[all]{xy}
\usepackage[mathscr]{eucal}
\usepackage{amsthm,amsmath,amssymb,amscd} 
\usepackage{epsfig}
\newtheorem{theo}{Theorem}

\newtheorem{ob}{Remark}

\begin{document}
\title{Relative equilibria in quasi-homogeneous three body problems}
\author{John A. Arredondo}
\date{}

\maketitle

\begin{abstract}
In this paper we find the families of relative equilibria for the three body problem in the plane, when the interaction between the bodies is given by a quasi-homogeneous potential, which is the sum of two homogeneous functions. The number of the relative equilibria depends of the values of the masses and of the size of the system, measured by the moment of inertia.
\end{abstract}

\section{Introduction}
In this paper we study a planar three body problem where the interaction between the bodies is given by a potential of the form
\begin{equation}\label{cuasi_pot}
U(r)=\frac{A}{r^\alpha}\pm\frac{B}{r^\beta},
\end{equation}
where $r$ is the distance between the bodies, $A$, $B$, $\alpha$ and $\beta$ are positive constants. This kind of potentials are called quasi-homogeneous because they are the sum of two functions which are homogeneous, in this case with homogeneity degree  $-\alpha$ and $-\beta$. Expression (\ref{cuasi_pot}) generalize many very well known quasi-homogeneous potentials as Birkhoff, Manev, Van der Waals, Libhoff, Schwarzschild,  Lennard-Jones, of course the classical  Newton and Coulomb and in some cases, potentials that come from exact solutions of the general relativity equations. In what follows our main porpoise is give a characterization of  the special periodic solutions called relative equilibria, associated to the famous problem of central configurations. Our main contribution is give an analytical proof of the kind of relative equilibria in two situations, the attractive-attractive case and the attractive-repulsive case, which for us means that in expression (\ref{cuasi_pot}) the components of the potential are both positive or one positive and other negative, respectively. Specifically we show that relative equilibria can correspond to arrangements of the bodies in equilateral, isosceles and scalene triangles, in function of the different values of the masses. 

This problem has been studied before in specific context by several authors, some introductory aspects can be found in \cite{C} for relative equilibria with Lennard-Jones potential in the two and three body problem with  equal masses. In our previous article \cite{AP} we study relative equilibria in the three body problem with Schwarszchild potential and any masses. This two previous references using numerical tolls in their conclusions.  In \cite{E} the authors give a proof of Moulton theorem for quasi-homogeneous potentials in general, and  in \cite{J, P} the authors explore the nature of the central configurations and their relation with the orbits of the bodies.

The article is organized as follows: In section 2 we introduce the equations of motion and the definitions relating to relative equilibria and central configurations, and along section 3 we study the planar relative equilibria for two cases: the attractive-repulsive where we show how the number of relative equilibria depends of the size of the system and the attractive-attractive, in both for all the different values of the masses. 

\section{Equations of motion}

We consider systems of  three bodies with masses $m_1, m_2, m_3$, moving in the $2$--dimensional Euclidean space under the influence of a quasihomoheneous type-potential. Let $\mathbf{q}_i\in \mathbb{R}^2$ denote the position of the particle $i$ in an inertial coordinate system and let $\mathbf{q}=(\mathbf{q}_1,\mathbf{q}_2,\mathbf{q}_3)$ the position vector, then the  generalized quasihomogeneous potential for the three body problem takes the form
\begin{equation}\label{potential3b}
U(\mathbf{q}) = \sum_{i \neq j}^3 \frac{A_{(m_im_j)}}{r_{ij}^\alpha} \pm  \sum_{i \neq j}^3 \frac{B_{(m_im_j)}}{r_{ij}^\beta},
\end{equation}
where $r_{ij}= | \mathbf{q}_i - \mathbf{q}_j |$, $\alpha$ and $\beta$ are positive constants for which we consider $\alpha>\beta$, and each $A_{(m_im_j)}$, $B_{(m_im_j)}$ is a positive constant depending of the interactions between the masses $m_i$ and $m_j$, respectively with $(i,j,k)$ permuting cyclically  in $(1,2,3)$, that is $(i,j,k) \sim (1,2,3)$. The equations of motion associated to the  potential
(\ref{potential3b}) are given by
\begin{equation}
\ddot{\mathbf{q}}=-\nabla U(\mathbf{q})\ , \label{eq1}
\end{equation}
Along the paper we assume as is classically, that the center of mass of the three particles is fixed at the origin, i.e.
\begin{equation}\label{cm}
\sum_{i=1}^3 m_i\mathbf{q}_i=0.
\end{equation}

The goal in this paper is the analysis of  the {\em relative equilibrium}; that is, solutions of (\ref{eq1}) that become equilibrium points in an uniformly rotating coordinate system (see (\cite{ME} for more details). Relative equilibria are characterized as follows: Let $R(\omega t)$ denote the $6\times 6$ block diagonal matrix with $3$ blocks of size $2\times 2$ corresponding to the canonical rotation in the plane. Let $\mathbf{x}\in (\mathbb{R}^6)$ be a configuration of the $3$ particles, and let $\mathbf{q}(t)=R(\omega t) \mathbf{x}$, where the constant $\omega$ is the angular velocity of the uniform rotating coordinate system. In the coordinate system $\mathbf{x}$ the equation of motion (\ref{eq1}) becomes
\begin{equation}
\label{eq2} \ddot{\mathbf{x}}+2\omega J \dot{\mathbf{x}}=-\nabla
U(\mathbf{x})+\omega^2 \mathbf{x}\ ,
\end{equation}
where $J$ is the $6\times 6$ block usual symplectic  matrix. A configuration $\mathbf{x}$ is called {\em central configuration} for system (\ref{eq1}) if and only if $\mathbf{x}$ is an equilibrium point of system (\ref{eq2}). That is, if and only if
 \begin{equation}\label{re}
 -\nabla U(\mathbf{x})+\omega^2 \mathbf{x}=\mathbf{0}\ ,
\end{equation}
for some $\omega$. If $\mathbf{x}$ is a central configuration, then 
\begin{equation}
\mathbf{q}(t)=R(\omega t)\mathbf{x}
\end{equation}
 is a  relative equilibrium solution of system (\ref{eq1}), which is also periodic with period $T=2\pi/|\omega|$. So when we obtain a central configuration, we are also getting the corresponding relative equilibria. Because of the explanation above, is very usual  talk about central configurations and relative equilibria as equivalent concepts.

For non-expert readers in this topic, is useful to remark that equation (\ref{re}) for a central configuration $\mathbf{q}=\mathbf{x}$, says that a central configuration in the space $\mathbf{q}$ is a configuration of the particles for which the particle $\mathbf{q}$ and the acceleration $\ddot{\mathbf{q}}$ vectors of each particle are proportional, with the same constant of proportionality $\omega^2$.

\section{Non-collinear central configurations}

Before star with the presentation of our results, is useful remember that if $u=f(x)$,  $x=(x_1,...x_n)$,$x_1=g_1(y),...,x_n=g_n(y)$ with $y=(y_1,...y_m)$,  $m\geq n$, and if $rank(A)=n$, where
                $$ A = \left(
                         \begin{array}{ccc}
                           \frac{\partial x_1}{\partial y_1} & \cdots & \frac{\partial x_n}{\partial y_n} \\
                          \vdots & \ddots &   \\
                           \frac{\partial x_1}{\partial y_m} & & \frac{\partial x_n}{\partial y_m} , \\

                         \end{array}
                       \right),$$
then $\nabla f(x)=0$ if and only if $\nabla u(y)=0$. This fact was previously used in  \cite{C}, and then for many authors for a similar porpoise.

\subsection{Attractive-repulsive case}

In this section we consider that the interaction between the bodies correspond to a quasi-homogeneous potential,  with one attractive component and other repulsive. So expression (\ref{potential3b}) is rewritten as
\begin{equation}\label{potential3bI}
U(\mathbf{q}) = \sum_{i \neq j}^3 \frac{A_{(m_im_j)}}{r_{ij}^\alpha} -  \sum_{i \neq j}^3 \frac{B_{(m_im_j)}}{r_{ij}^\beta}.
\end{equation}
Because at this case, central configurations are not invariant under homotheties,  is natural to think that the number of these depend of the size of the system measured by the moment of inertia $I$, which  can be written in terms of the mutual distances as
\begin{equation}\label{moment of I}
I=\frac{1}{M}(m_1m_2r_{12}^{2}+m_1m_3r_{13}^{2}+m_2m_3r_{23}^{2}),
\end{equation}
where  $M=m_1+m_2+m_3$  and we have assumed, without loss of generality, that the value of the three equal masses is $\dfrac{1}{3}$ and then $M=1$. The main result in this case is the following.

 \begin{theo}\label{theo1}
Consider the planar $3$--body problem, where the mutual interaction between the particles is given by a quasi-homogeneous potential (\ref{potential3bI}), then,
\begin{enumerate}
\item If the three masses are equal, the relative equilibria can be equilateral or isosceles triangles. The number of relative equilibria depend of the moment of inertia, and there are four bifurcation values for $I.$
\item  If two  masses are equal, the relative equilibria  correspond to isosceles or scalene triangles.
\item If the three masses are different, the relative equilibria  are always scalene triangles
\end{enumerate}
 \end{theo}

{\bf Proof:}  Central configurations are solutions of the system
\begin{eqnarray}{}\label{planarcc1}
-\frac{\alpha A_3}{r_{12}^{\alpha+1}}+\frac{\beta B_3}{r_{12}^{\beta+1}}+r_{12}\omega^2 &=&0, \nonumber\\
-\frac{\alpha  A_2}{r_{13}^{\alpha+1}}+\frac{\beta B_2}{r_{13}^{\beta+1}}+r_{13}\omega^2 &=&0,\\
-\frac{\alpha  A_1}{r_{23}^{\alpha+1}}+\frac{\beta B_1}{r_{23}^{\beta+1}}+r_{23}\omega^2 &=&0,\nonumber\\
\frac{1}{9}(r_{12}^{2}+r_{13}^{2}+r_{23}^{2}) &=&I,\nonumber
\end{eqnarray}
the first three equations are equivalent to
\begin{equation}\label{omega_cuasi}
\frac{\alpha A_3}{r_{12}^{\alpha+2}}-\frac{\beta B_3}{r_{12}^{\beta+2}}=\frac{\alpha A_3}{r_{13}^{\alpha+2}}-\frac{\beta B_3}{r_{13}^{\beta+2}}=
\frac{\alpha A_3}{r_{23}^{\alpha+2}}-\frac{\beta B_3}{r_{23}^{\beta+2}}=\omega^2,
\end{equation}
which for our convenience, abusing of notation we rewrite as
\begin{equation}\label{omega_cuasi1}
\frac{A_3}{r_{12}^{\alpha+2}}-\frac{B_3}{r_{12}^{\beta+2}}=\frac{A_2}{r_{13}^{\alpha+2}}-\frac{B_2}{r_{13}^{\beta+2}}=
\frac{A_1}{r_{23}^{\alpha+2}}-\frac{ B_1}{r_{23}^{\beta+2}}=\omega^2.
\end{equation}
With the introduction above, let us star with the proof of the first part:

1.  Let $f(x)=\dfrac{A}{x^{\alpha+2}}-\dfrac{B}{x^{\beta+2}}$ be the equivalent of equations (\ref{omega_cuasi1}) for three equal masses. Since $x$ represents a distance, it is enough to analyze the function $f(x)$ for positive $x$. At this case the principal remarks of $f(x)$ are:
\begin{figure}[h]
\centerline
{\includegraphics[scale=0.7]{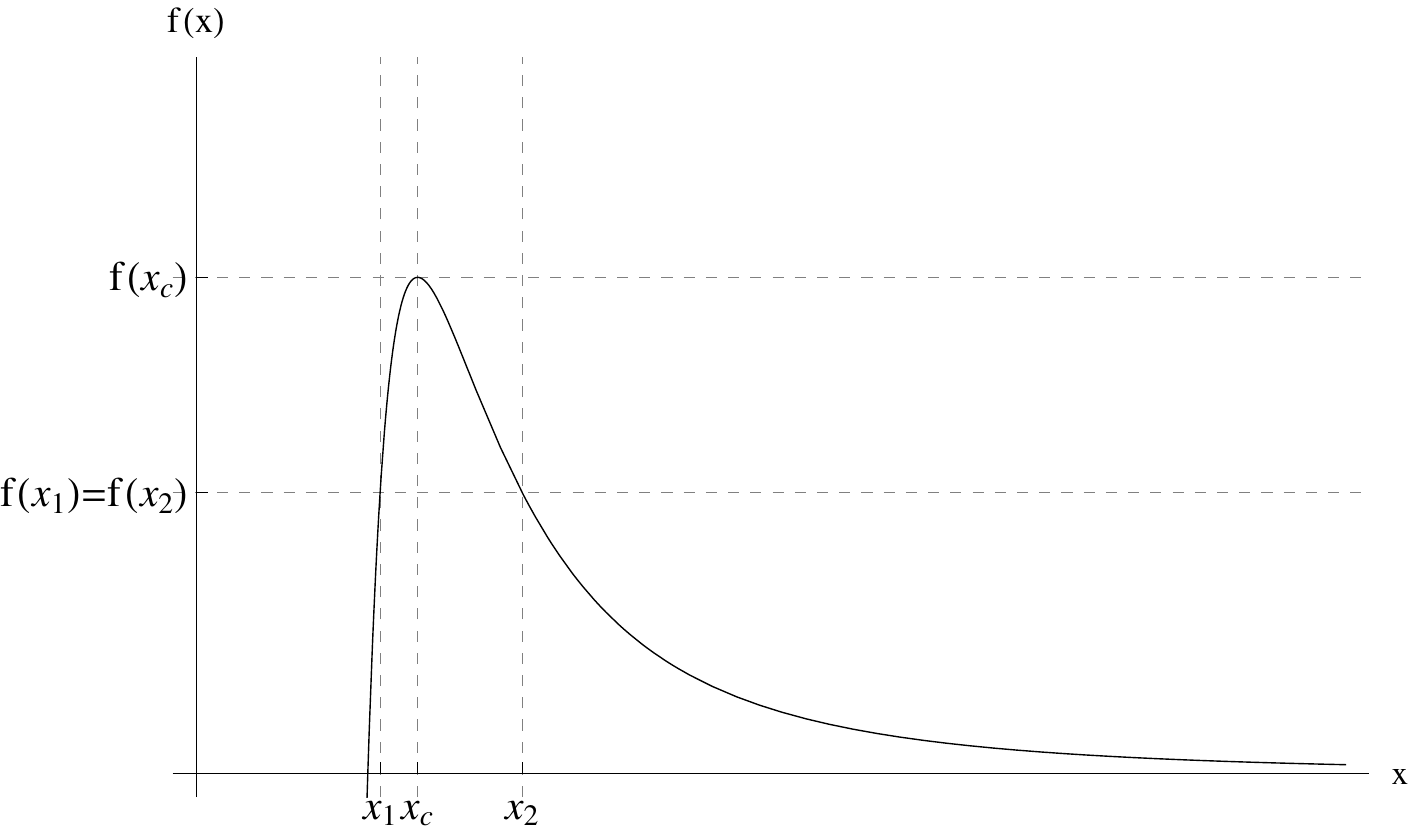}}
\caption{ \label{B_neg} \emph{The graph of $f(x)$, with $A, \alpha, \beta$ and $B$  positive constants.}}
\end{figure}
\begin{enumerate}
\item[1)] the function is zero at $x_{0}=(B/A)^\frac{1}{\beta - \alpha}$, 
\item[2)] has a maximun point at $x_{c}=((\beta+2) B/(\alpha+2) A)^\frac{1}{\beta - \alpha}$,
\item[3)] the limites are $\lim_{_{x \to0}}f(x)=-\infty$ y  $\lim_{_{x \to \infty}}f(x)=0$.
\end{enumerate}
Therefore, a generic graph of the function $f(x)\geq 0$, is presented in figure \ref{B_neg}, with this, the study of equations (\ref{omega_cuasi1}) can be reduced to analyze 
\begin{equation}\label{same}
f(r_{12})=f(r_{13})=f(r_{23})=\omega^2.
\end{equation}
Fixed a value of $\eta \in (0,\beta)$ where $\beta=f(x_c)$ we can find two different values, $x_1 \in (x_0,x_c)$ and $x_2 \in (x_c,+\infty)$, satisfying $f(x_1)=f(x_2)= \eta$. So with these two possible solutions, we can get equilateral solutions, when the configuration of the three particles has equal sides $r_{12}=r_{13}=r_{23}=x_1$ or $r_{12}=r_{13}=r_{23}=x_2$, and isosceles solutions when the configuration of the three particles has  two sides $r_{ij}=x_1$ and the third one is equal to $x_2$ or when two sides are $r_{ij}=x_2$ and the third one is equal to $x_1$ . When $\eta = \beta$, there is a unique value $x_c$ such that $f(x_c)=\beta$, only in this case, we have a unique solution of (\ref{same}) given by $r_{12}=r_{13}=r_{23}=x_{c}$.
\begin{figure}[h]
\centerline
{\includegraphics[scale=0.8]{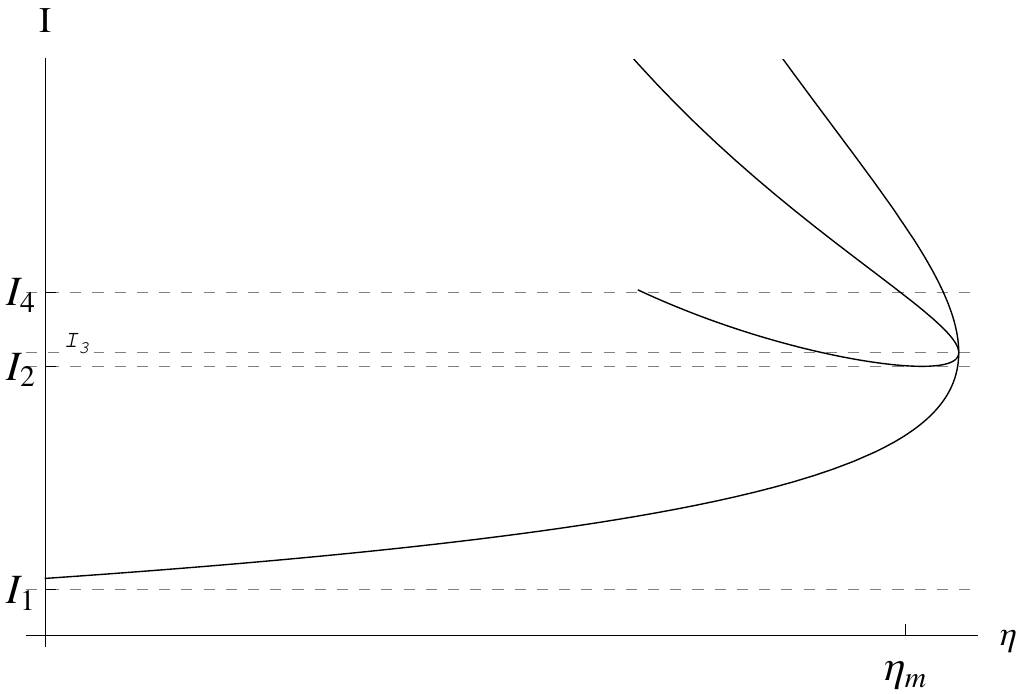}}
\caption{ \label{iner_moment} \emph{Behavior of the moment of inertia for $f(x)=\dfrac{A}{x^{\alpha}}-\dfrac{B}{x^{\beta}}$ with $\alpha > \beta$. The ascendant curve corresponds to equilateral solutions of type (1), the curve between $I_2$ and $I_4$ corresponds to isosceles solutions of type (4), the next curve corresponds to isosceles solutions of type (3) and the upper curve corresponds to equilateral solutions of type (2).}}
\end{figure}

To prove the second part of the statement one, we have to analyze the moment of inertia for the above solutions and check that all solutions satisfy the triangular inequality. For equilateral solutions, the moment of inertia is given by
\begin{itemize}
\item[1)]  $I=\dfrac{1}{3}x_{1}^2$, \quad \text{where} \quad $x_1\in (x_0, x_c]$,
\item[2)]  $I=\dfrac{1}{3}x_{2}^2$,  \quad \text{where} \quad $x_2\in [x_c, \infty)$,
\end{itemize}
and for isosceles solutions we get
\begin{itemize}
\item[3)]  $I=\dfrac{1}{9}(2x_{2}^2+x_1^2)$, \quad \text{where always} \quad $2x_i>x_j,$\, $i,j=1,2$, $i\neq j$,
\item[4)]  $I=\dfrac{1}{9}(2x_{1}^2+x_2^2)$, \quad \text{where at some value} \quad $2x_1=x_2$.
\end{itemize}
In figure \ref{iner_moment}, we have depicted these four families of solutions. For values $I \leq I_1$, there are not relative equilibria, since for those  the function $f(x)$ is negative. A simple inspection verify that in the first three expressions, the distance $x_1$ and $x_2$ fulfill the triangular inequality. In the four one, because  $x_1$ is restricted to the set $(x_0, x_c]$ and $x_2$ can take values in the set $ [x_c, \infty)$, at some point will happen that $x_2\geq 2x_1$, therefore, at this point these distances do not satisfy the triangular inequality. We want estimate when $x_2=2x_1$, for this, let be $x_1=x_c-\epsilon$ and $x_2=x_c+\epsilon$, evaluate the inertia moment at these values and compare with the inertia moment evaluated at $x=x_c$. With this we get
\begin{equation}\label{eq_11}
I=3x_c^2-2\epsilon x_c+3\epsilon^2 \leq 3x_c^2,
\end{equation}
expression (\ref{eq_11}) is valid if $\epsilon \leq \frac{2}{3}x_c$, i.e., for $\epsilon$ small, the inertia moment  given by  $9I=2x_{1}^2+x_2^2$, is a decreasing function. If $\epsilon \geq \frac{2}{3}x_c$, then $I\geq 3x_c$ and therefore at some point the curve change from decreasing to increasing. To determine if this happens, let us make an idea or the location of $\frac{2}{3}x_c$ comparing it, with the point where $f(x)$ is zero, i.e., we want to establish if $\frac{2}{3}x_c>x_0$ or $\frac{2}{3}x_c<x_0$. After some calculations, we find that
\begin{equation}
\left(\frac{\beta+2}{\alpha+2}\right)<\left(\frac{3}{2}\right)^{\beta - \alpha}
\end{equation}
for all $\alpha\geq 1$, $\beta>\alpha$. Therefore  $\frac{2}{3}x_c<x_0$, which imply that $x_c-x_0\leq \frac{1}{3}$. 

Taking values of $x_1$ and $x_2$ in the ball with center in $x_c$ and radio $\frac{1}{2}x_c$ $\mathcal{B}(x_c, \frac{1}{2}x_c)$, we find that for values at the left in the limit case $x_1=x_0$, whit this and taking $x_2= \frac{3}{2}x_c$, the triangular inequality is always fulfill, this means that
\begin{equation}
2\left(\frac{B}{A}\right)^\frac{1}{\beta - \alpha}\geq \frac{3}{2}\left(\frac{(\beta+2)B}{(\alpha+2)A}\right)^\frac{1}{\beta - \alpha},
\end{equation}
since this expression can be simplified as
$$\left(\frac{4}{3}\right)^{\beta - \alpha}\geq\left(\frac{\beta+2}{\alpha+2}\right)^{\beta - \alpha},$$
which is valid for all $\alpha\geq 1$, $\beta-\alpha\geq 1$. Hence, if is true for $x_1=x_0$, also is valid for a different $x_1>x_0$. Now we want to establish if in the ball $\mathcal{B}(x_c, \frac{1}{2}x_c)$ there are $x_1$ and $x_2$ such that, the inertia moment evaluated at those values, is greater than the inertia moment evaluated at $x_c$. For this we take the limit values for the right, which is $x_2=\frac{3}{2}x_c$, and we looking for a $x_1$, which verify the inequality
\begin{equation}
2x_1^2+\frac{9}{4}x_c^2\geq 3x_c^2,
\end{equation}
this imply that $x_1\geq \frac{\sqrt{6}}{4}x_c$, but $x_0\geq \frac{2}{3}x_c> \frac{\sqrt{6}}{4}x_c$. Since this is true for $x_0$, also is true for all $x_1\geq x_0$. Then for values in $\mathcal{B}(x_c,\frac{1}{2}x_c)$ and for a certain number of values outside of this ball, always can be found $x_1$ and $x_2$ which satisfy
\begin{equation}
2x_1^2+ x_2^2\geq 3x_c^2,
\end{equation}
therefore, the function $9I=2x_{1}^2+x_2^2$ at some point changes from decreasing to increasing, i.e., the function has a critical point, which is a minimum, whereby for any positive constants $A$, $B$, $\beta>\alpha\geq 1$ the figure \ref{iner_moment} is generic. \qed

\medskip

Now we proceed with the proof of the second part of the theorem for the case of two equal masses. 

\begin{figure}[h]\vspace{-11cm}\hspace{2cm}
\includegraphics[width=22cm, bb=0 0 640 480]{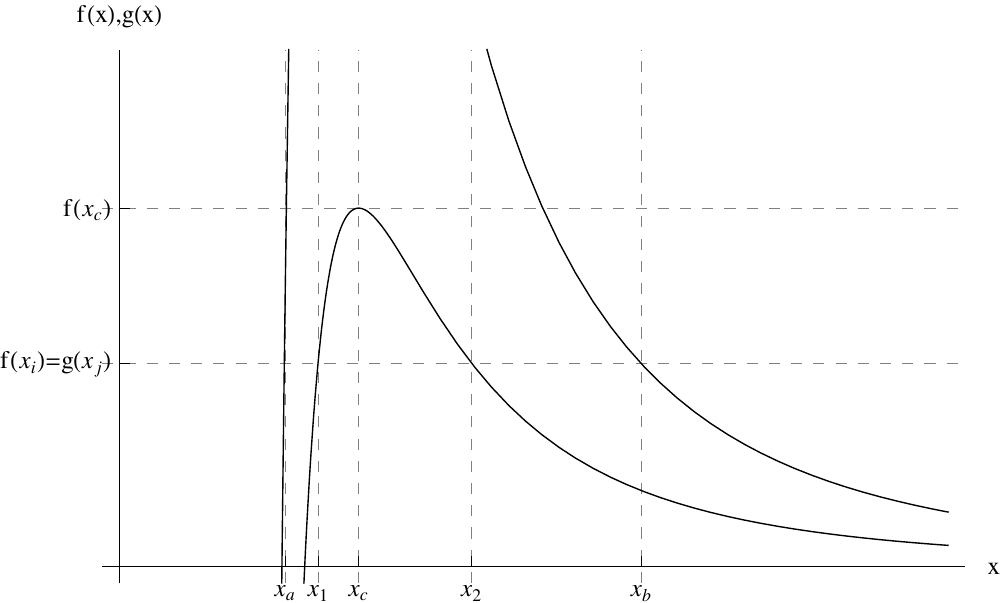}
\caption{ \label{Graf. 4} \emph{Graph of $f(x)$ and $g(x)$ for two equal masses}}
\end{figure}

2.  Equations (\ref{same}) when two masses are equal becomes:

\begin{align}\label{factork}
\frac{A}{r_{12}^{\alpha+2}}-\frac{B}{r_{12}^{\beta+2}}=\frac{A}{r_{13}^{\alpha+2}}-\frac{B}{r_{13}^{\beta+2}}=
\frac{kA}{r_{23}^{\alpha+2}}-\frac{k B}{r_{23}^{\beta+2}}=\omega^2,
\end{align}
where $k$ is a positive constant that measures how the third mass is increasing with respect to the other two, the analysis when the third mass is decreasing is similar. We denote the first two identical functions in the above equation as $f(x)$ and the third one is represented as $g(x)$, $(g(x)=kf(x))$, see Fig. \ref{Graf. 4}.

Both functions $f(x)$ and $g(x)$ have a maximum at $x_c$. Let $\beta=f(x_c)$, for any $\eta \in (0,\beta]$, the straight line $y=\eta$ intersects the graph of $f(x)$ in two points that we call $f(x_1)=f(x_2)=\eta$ and the graph of $g(x)$ in the points $g(x_a)=g(x_b)=\eta$ see Fig \ref{Graf. 4}. Then for this value of $\eta$ we can have isosceles configurations if $r_{12}=r_{13}=x_1, r_{23}=x_a$, $r_{12}=r_{13}=x_1, r_{23}=x_b$, $r_{12}=r_{13}=x_2, r_{23}=x_a$ and $r_{12}=r_{13}=x_2, r_{23}=x_b$ or scalene configurations if $r_{12}=x_1, r_{13}=x_2, r_{23}=x_a$ and $r_{12}=x_2, r_{13}=x_1, r_{23}=x_b$.

\begin{figure}[h]\vspace{-11cm}\hspace{3cm}
\includegraphics[width=21cm, bb=0 0 640 480]{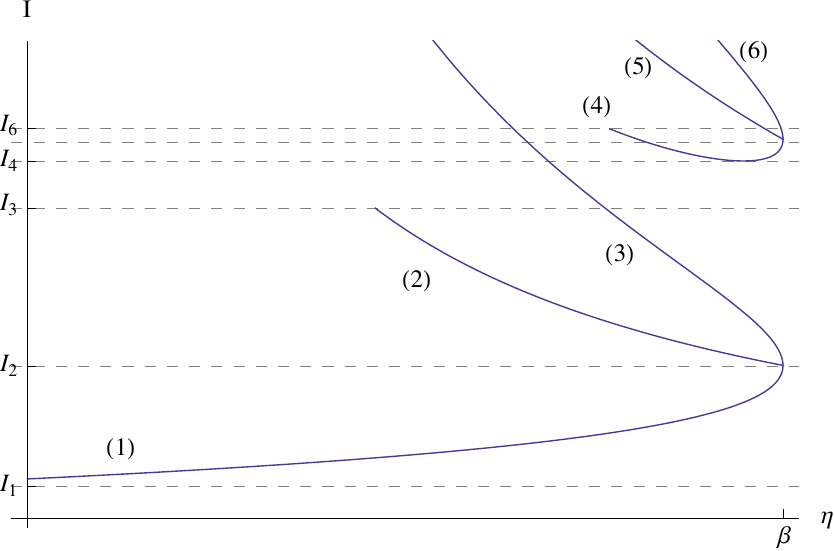}
\caption{ \label{Graf.5} \emph{Behavior of the moment of inertia for two equal masses.}}
\end{figure}

The next step is to analyze the moment of inertia for the solutions  above. Since the constants $A$ and $B$ are arbitrary in our analysis, we can include the value of the masses in these parameters, so when we take different masses, we can modify the value of the parameters and do the study  considering in all cases that the particles have equal masses. The curves that we are representing correspond to:
\begin{itemize}
\item (1) $9I=2x_{1}^{2}+x_{a}^{2}$, for the interval $(I_1,I_2]$.
\item (2) $9I=x_{1}^{2}+x_{2}^{2}+x_{a}^{2}$, for the interval $[I_2,I_3]$.
\item (3) $9I=2x_{2}^{2}+x_{a}^{2}$, for the interval $(I_2,\infty)$.
\item (4) $9I=2x_{1}^{2}+x_{b}^{2}$, for the interval $[I_4,I_6]$.
\item (5) $9I=x_{1}^{2}+x_{2}^{2}+x_{b}^{2}$, for the interval $[I_5,I_7]$.
\item (6) $9I=2x_{2}^{2}+x_{b}^{2}$, for the interval $(I_5,\infty)$.
\end{itemize}
Again we remark that all solutions must satisfy the triangular inequality. For expressions $9I=2x_{1}^{2}+x_{a}^{2}$ and  $9I=2x_{2}^{2}+x_{a}^{2}$ this condition is always fulfill, because $x_c-x_0\leq \frac{1}{3}$. In each of the other cases we have to analyze that condition. In figure \ref{Graf.5} we represent all this curves. 

For $9I=x_{1}^{2}+x_{2}^{2}+x_a^2$ the triangular inequality fail if $x_a+x_1\leq x_2$. At the initial point $x_1=x_2=x_c$ and this function take the value $9I=2x_c^2+x_a^2$. If we move each of this values a quantity $\epsilon$, $x_1=x_c-\epsilon$, $\tilde{x}_a=x_a-\epsilon$ and $x_2=x_c+\epsilon$ and evaluating the function we have $9I=2x_c^2+x_a^2+\epsilon(3\epsilon-2x_a)$, therefore, close to the initial point, this function is decreasing, gets a minimum point and then is increasing for all time.

For $9I=2x_{1}^{2}+x_{b}^{2}$ the triangular inequality fail if $x_b\geq 2x_1$, actually, this curve is possible  only if $x_b< 2x_c$. As before, close to the initial values, $x_1=x_c-\epsilon$ and $\tilde{x}_b=x_b-\epsilon$, evaluating the function we have $9I=2x_c^2+x_b^2+\epsilon(3\epsilon+2x_b-4x_c)$, therefore, close to the initial point, this curve decreases, gets a minimum point and then  increase monotonically, until $x_b = 2x_1$. the value  where this curves ends is always large than the values where it stars. 

For $9I=x_{1}^{2}+x_{2}^{2}+x_b^2$ the triangular inequality fail if $x_2+x_1\leq x_b$. At the initial point $x_1=x_2=x_c$ and this function take the value $9I=2x_c^2+x_b^2$. If we move each of this values a quantity $\epsilon$, $x_1=x_c-\epsilon$, $\tilde{x}_b=x_b-\epsilon$ and $x_2=x_c+\epsilon$ and evaluating the function we have $9I=2x_c^2+x_b^2+\epsilon(3\epsilon+2x_b)$, which imply that $2x_c^2+x_b^2+\epsilon(3\epsilon+2x_b)>2x_c^2+x_b^2$, therefore this function is monotonically increasing.

For $9I=2x_{2}^{2}+x_{b}^{2}$ the triangular inequality fail if $x_b\geq 2x_2$, analyzing close to the initial values as before, we find that the curve is monotonically increasing.

\begin{ob}
To establish a estimation of   the "'distance'' between $f(x)$ and $g(x)=kf(x)$ for which all the curves study have sense, let as analyze the extreme case where $x_b=2x_c$. With this values in the expressions of $f(x)$ and $g(x)$, we get that $k$ denoted $\tilde{k}$ for this special values, in terms of the other parameters, is given by
\begin{equation}
\tilde{k}=2^{\beta+2}\frac{Ax_c^{\beta-\alpha}-B}{2^{\beta-\alpha}Ax_c-B}
\end{equation}
therefore if $k<\tilde{k}$, then all the six possible expressions for the inertia moment, analyzed before, exist for some interval in which they fulfill the triangular inequality.  If $k<\tilde{k}$, then the expressions for the inertia moment
\begin{itemize}
\item (4) $9I=2x_{1}^{2}+x_{b}^{2}$, 
\item (5) $9I=x_{1}^{2}+x_{2}^{2}+x_{b}^{2}$, 
\item (6) $9I=2x_{2}^{2}+x_{b}^{2}$, 
\end{itemize}
do not fulfill the triangle inequality and therefore they can not appear in figure \ref{Graf.5}. So the number of central configurations for two equal masses can be change from a minimum of three families, two isosceles and one scalene, until a maximum of six families, four isosceles and two scalene.
\end{ob}
\qed

\medskip

Now we proceed with the proof of the last part of the theorem for the case of three different masses. 

\begin{figure}[h]\vspace{-5.5cm}\hspace{2cm}
\includegraphics[width=15cm, bb=0 0 640 480]{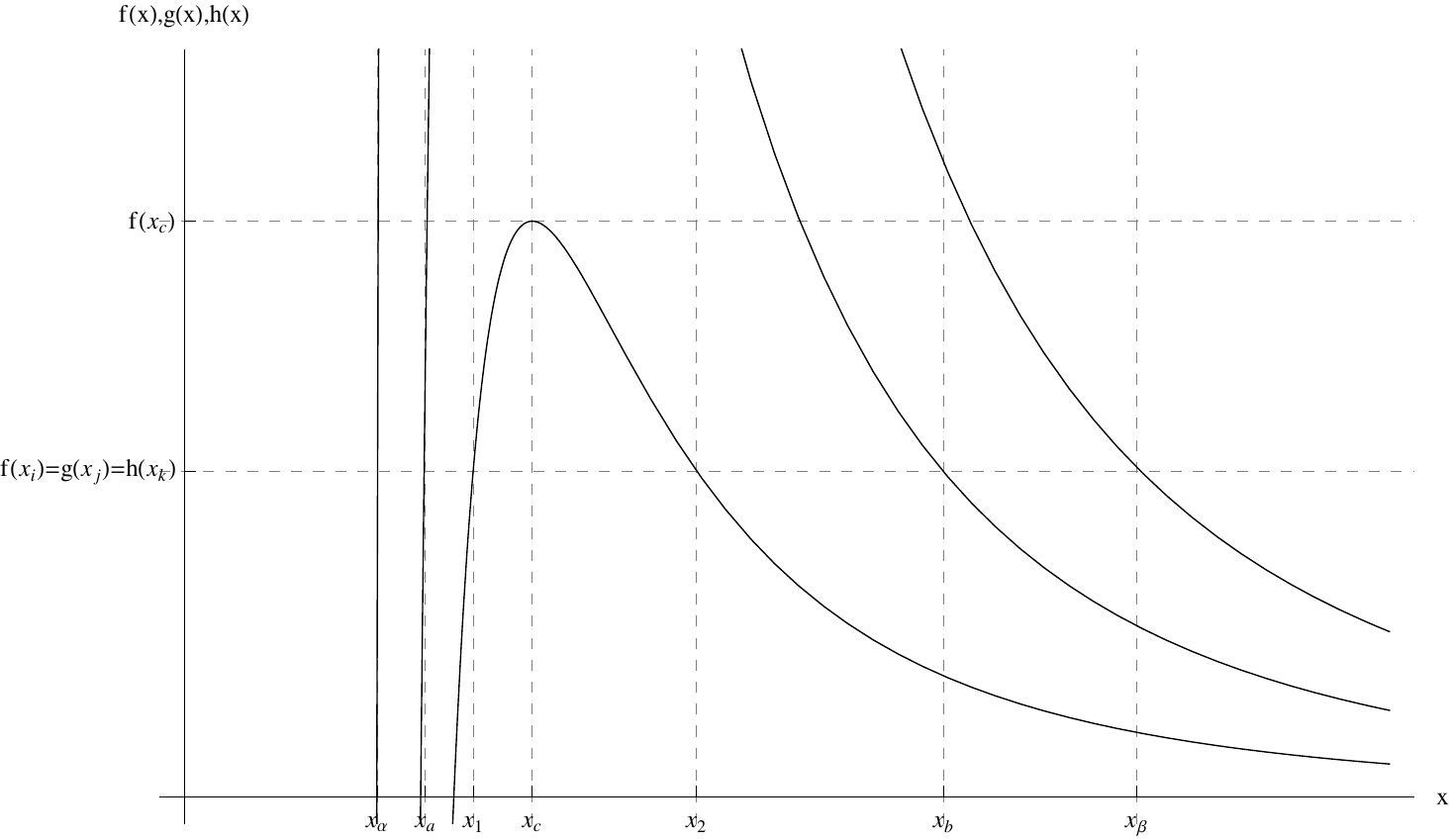}
\caption{ \label{Graf. 7} \emph{The graph of $f(x), g(x)=kf(x)$ and $h(x)=k_1f(x)$.}}
\end{figure}

3.  Equations (\ref{same}) when the three masses are different becomes:
\begin{align}\label{factork1}
\frac{A}{r_{12}^{\alpha+2}}-\frac{B}{r_{12}^{\beta+2}}=\frac{kA}{r_{13}^{\alpha+2}}-\frac{kB}{r_{13}^{\beta+2}}=
\frac{k_1A}{r_{23}^{\alpha+2}}-\frac{k_2 B}{r_{23}^{\beta+2}}=\omega^2,
\end{align}
or equivalent
\begin{align}
f(r_{12})=k f(r_{13})=k_1 f(r_{23})=\omega^2.
\end{align}
 with $k$ and $k_1$ positive constants. By using this notation $g(x)=kf(x)$ and $h(x)=k_1f(x)$. Here the graphs of the three functions are different, but all of them reach the maximum at the same point $x_c$. In Fig.  \ref{Graf. 7},  we have plotted the graphs of these functions.

\begin{figure}[h] \vspace{-11.5cm}\hspace{3cm}
\includegraphics[width=22cm, bb=0 0 640 480]{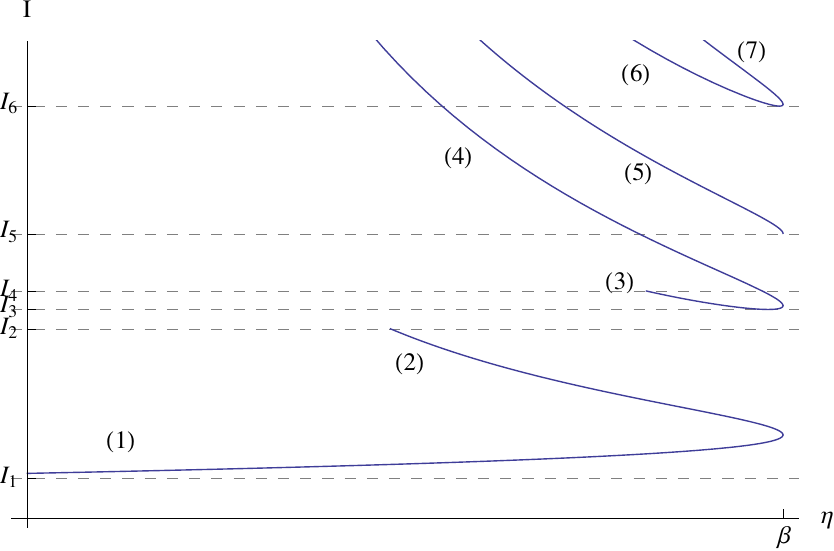}
\caption{ \label{Dib.4} \emph{Behavior of moment of inertia in the case that the three masses are different}}
\end{figure}

As in the previous two cases, the straight line $y=\eta$ intersects every graph in two points that we call $f(x_1)=f(x_2)=\eta$, $g(x_a)=g(x_b)=\eta$ and $h(x_\alpha)=h(x_\beta)=\eta$. Then for this value of $\eta$, in our problem of find central configurations, all solutions correspond to scalene triangles. The corresponding curves, associated to the inertia moment and the corresponding interval of definition (see figure \ref{Dib.4}) are

\begin{itemize}
\item (1) $I=\frac{1}{9}(x_{\alpha}^{2}+x_{a}^{2}+x_{1}^{2})$ in the interval $(I_1,I_2),$
\item (2) $I=\frac{1}{9}(x_{\alpha}^{2}+x_{a}^{2}+x_{2}^{2})$ in the interval $(I_1,I_2],$
\item (3) $I=\frac{1}{9}(x_{\alpha}^{2}+x_{b}^{2}+x_{1}^{2})$ in the interval $[I_3,I_4],$
\item (4) $I=\frac{1}{9}(x_{\alpha}^{2}+x_{b}^{2}+x_{2}^{2})$ in the interval $(I_4,I_8],$
\item (5) $I=\frac{1}{9}(x_{\beta}^{2}+x_{a}^{2}+x_{2}^{2})$ in the interval $[I_5,I_7],$
\item (6) $I=\frac{1}{9}(x_{\beta}^{2}+x_{b}^{2}+x_{1}^{2})$ in the interval $[I_6,I_9]$,
\item (7) $I=\frac{1}{9}(x_{\beta}^{2}+x_{b}^{2}+x_{2}^{2})$, in the interval $(I_6,\infty)$.
\item (7) $I=\frac{1}{9}(x_{\beta}^{2}+x_{a}^{2}+x_{1}^{2})$, in the interval $[I_5,I_6]$,.
\end{itemize}

Verifying the triangular inequality for the expressions above, and taking into account that $x_1\in [\frac{2}{3}x_c,x_c]$, $x_2\in [x_c, \infty)$, $x_a\in [\frac{2}{3}x_c,x_c)$, $x_b\in (x_c, \infty$, $x_\alpha \in [\frac{2}{3}x_c,x_c)$ and $x_\beta\in (x_c,\infty)$, with $x_\alpha<x_a<x_1$ and $x_2<x_b<x_\beta$, we find the following resolves: 

For $I=\frac{1}{9}(x_{\alpha}^{2}+x_{a}^{2}+x_{1}^{2})$ the triangular inequality failure if $x_\alpha +x_a<x_1$, and because the largest value of $x_1$ is $x_c$ and the smaller for the other two values are $x_c/3$, it follows that $4/3x_c<x_1=x_c$, which is a contradiction, so the values in this curve always fulfill the triangular inequality.   (see  Fig. \ref{Dib.4}).

For $I=\frac{1}{9}(x_{\alpha}^{2}+x_{a}^{2}+x_{2}^{2})$  the triangular inequality failure if $x_\alpha +x_a<x_1$, but because the minimum values $min x_\alpha=min x_a=\frac{2}{3}x_c$ and $min x_2=x_c$, there is always a small interval in which this curve is well define (i.e., the inertia moment curve take values which fulfill the triangular inequality), and ends in a bifurcation value of collinear configuration.

For $I=\frac{1}{9}(x_{\alpha}^{2}+x_{b}^{2}+x_{1}^{2})$  the triangular inequality failure if $x_\alpha +x_1<x_b$, in this case with  $min x_\alpha=min x_1=\frac{2}{3}x_c$, we find that if $x_b<\frac{4}{3}x_c$ there is a small  interval in which this curve is well define and ends in a bifurcation value of collinear configuration. If  $x_b>\frac{4}{3}x_c$ the triangular inequality is not fulfill, and the curve does not exist.

For $I=\frac{1}{9}(x_{\alpha}^{2}+x_{b}^{2}+x_{2}^{2})$ the triangular inequality failure if $x_\alpha +x_2<x_b$, we find that if $x_b<\frac{5}{3}x_c$ when $x_2=x_c$, the corresponding curve is well define. If  $x_b>\frac{5}{3}x_c$ the triangular inequality is not fulfill, and the curve does not exist.

For $I=\frac{1}{9}(x_{\beta}^{2}+x_{a}^{2}+x_{2}^{2})$  the triangular inequality failure if $x_a+x_2<x_\beta$, we find that if $x_\beta<\frac{5}{3}x_c$ when $x_2=x_c$, the corresponding curve is well define. If  $x_\beta>\frac{5}{3}x_c$ the triangular inequality is not fulfill, and the curve does not exist.

For $I=\frac{1}{9}(x_{\beta}^{2}+x_{b}^{2}+x_{1}^{2})$  the triangular inequality failure if $x_1+x_b<x_\beta$, we find that if the difference $x_\beta-x_b\in (\frac{2}{3}x_c,x_c)$, the corresponding curve is well define. Otherwise, the triangular inequality is not fulfill, and the curve does not exist. 

For $I=\frac{1}{9}(x_{\beta}^{2}+x_{b}^{2}+x_{2}^{2})$ the triangular inequality failure if $x_2+x_b<x_\beta$, If  $x_\beta-x_b<x_c$ when $x_2=x_c$ the corresponding curve is well define. Otherwise, the triangular inequality is not fulfill, and the curve does not exist. 

For $I=\frac{1}{9}(x_{\beta}^{2}+x_{a}^{2}+x_{1}^{2})$  the triangular inequality failure if $x_a+x_1<x_\beta$,  If  $x_\beta<\frac{4}{3}x_c$, there is a small  interval in which this curve is well define. . Otherwise, the triangular inequality is not fulfill, and the curve does not exist. 

\begin{ob}
As in the case of two equal masses, there is a "'distance'' between $f(x)$, $g(x)=kf(x)$ and $h(x)=k_1f(x)$ for which all the eight curves of the inertia moment, found for three different masses have sense. In the extreme case where $x_\beta=2x_c$, the corresponding $k_1$ denoted $\tilde{k}_1$, in terms of the other parameters, is given by
\begin{equation}
\tilde{k}_1=2^{\beta+2}\frac{Ax_c^{\beta-\alpha}-B}{2^{\beta-\alpha}Ax_c-B}
\end{equation}
therefore if $k<k_1<\tilde{k}_1$, then all the eight possible expressions for the inertia moment, analyzed before, exist for some interval in which they fulfill the triangular inequality.  If $k_1>\tilde{k}_1$, then the expressions for the inertia moment $I=\frac{1}{9}(x_{\alpha}^{2}+x_{a}^{2}+x_{1}^{2})$ is the only which vales fulfill the triangle inequality and therefore the only curve which appear in figure \ref{Dib.4}. So the number of central configurations for three different masses can be change from a minimum of one family until a maximum of eight families, scalene in all cases.
\end{ob}
\qed

\subsection{Attractive-attractive case.}

Briefly we present  the planar central configuration or the relative equilibria for quasi-homogeneous three body problem, when both components of the potential are attractive. The main result is the following.
\begin{theo}\label{T4} Consider the planar 3-body problem, where the mutual interactions between the particles is given by a quasi-homogeneous potential of the form
\begin{equation}\label{potential3bII}
U(\mathbf{q}) = \sum_{i \neq j}^3 \frac{A_{(m_im_j)}}{r_{ij}^\alpha} +  \sum_{i \neq j}^3 \frac{B_{(m_im_j)}}{r_{ij}^\beta},
\end{equation}
with $A_{(m_im_j)}, B_{(m_im_j)}$ positives constants. In this case the relative equilibria must belong to one of the following families and there are not bifurcation values with respect to the moment of inertia:
\begin{itemize}
\item If the three masses are equal, any relative equilibria must be arrangements of the particles in equilateral triangles.
\item If two of the three masses are equal, any relative equilibria must be arrangements of the particles in isosceles triangles.
\item If the three masses are different, any relative equilibria must be arrangements of the particles in scalene triangles.
\end{itemize}
\end{theo}

{\bf Proof:} The statements above follows from solving equations
\begin{align}\label{Eq-positivas}
f(r_{12})=kf(r_{13})=k_1 f(r_{23})=\omega^2,
\end{align}
where we are thinking this expressions as $f(r_{12})=f(x)=\dfrac{A}{x^{\alpha+2}}+\dfrac{B}{x^{\alpha+2}}$, $kf(r_{13})=g(x)$ and $k_1f(r_{23})=h(x)$, with $k, k_1>0$ and using that  $\lim_{_{x \to 0}}f(x)=\infty$ and  $\lim_{_{x \to \infty}}f(x)=0$.
\qed

The statement of theorem \ref{T4} includes those potentials which have been proposed as corrections, in the classical or in the relativistic context, of the Newtonian potential. For example our previous work \cite{AP} about central configurations with the Schwarszchild potential or the same situation with the Manev potential. On the other hand also generalizes the study of relative equilibria in restricted problems where the geometry of the bodies is involved \cite{A}, and the potential used is also a quasi-homogeneous function.

\begin{ob}
Because the statements of theorems \ref{theo1} and  \ref{T4} are independent of the values of the positive constants $A$ and $B$, in both cases $B=0$ and $\alpha=1$ implies that we recover the Newtonian potential, so we conclude that any positive or negative, small perturbation on the Newtonian potential change the classical resolve about relative equilibria in the three body problem, where all the planar relative equilibria correspond to arrangements in equilateral triangles, and for quasi-homogeneous potentials, this resolve persist only for equal masses.
\end{ob}

\begin{ob}
To conclude the whole problem, we have remember that collinear relative equilibria, are determined by the well known Moulton theorem, which for quasi-homogeneos potentials in the attractive-attractive case is probed for n-bodies in \cite{E}, and sets that the number of relative equilibria on the line is $\dfrac{n!}{2}$. With minors changes in the proof found there, for the  attractive-repulsive case can be proved that the number of relative equilibria on the line is $n!$. 
\end{ob}
\section{Acknowledgements}
We thank to Dr. Jesus Mci\~no for give us the idea of this paper and for his helpful comments.

\end{document}